\documentclass[a4paper,twoside]{article}
\usepackage{amsmath,amssymb,amsfonts}
\usepackage{graphics}
\newtheorem{proposition}{Proposition}[section]

\def\R{{\mathbb R}}

\date{\small\em \today}

\begin{document}
\title{Neighbourhoods of independence for random processes}
\author{ Khadiga Arwini and C.T.J. Dodson \\
{\small Department of Mathematics}\\
{\small University of Manchester Institute of Science and Technology} \\
{\small Manchester M60 1QD, UK}   }

\maketitle

\begin{abstract}
The Freund family of distributions becomes a Riemannian
4-manifold with Fisher information as metric;
we derive the induced $\alpha$-geometry, i.e., the $\alpha$-curvature,
$\alpha$-Ricci curvature with its eigenvales and eigenvectors,
the $\alpha$-scalar curvature etc.  We
show that the Freund manifold has a positive constant $0$-scalar
curvature, so geometrically it constitutes part of a sphere.
We consider special cases as submanifolds and discuss their
geometrical structures; one submanifold yields examples of
neighbourhoods of the independent case for  bivariate
distributions having identical exponential marginals. Thus, since
exponential distributions complement Poisson point processes, we
obtain a means to discuss the neighbourhood of independence for
random processes.
\end{abstract}

\section{Differential geometry of the Freund 4-manifold}
In ~\cite{gamran} we proved that every neighbourhood of an exponential
distribution contains a neighbourhood of gamma distributions, in
the subspace topology of $\R^3$ using information geometry and the
affine immersion of Dodson and Matsuzoe~\cite{affimm}.  For
general references on information geometry, see Amari et
al.~\cite{Am},~\cite{AmNag}. As part of a study of the information
geometry and topology of gamma and bivariate stochastic processes
cf. e.g.~\cite{amino},~\cite{ijtp},~\cite{gsis}, we have calculated the
geometry of the family of Freund bivariate mixture exponential
density functions. The importance of this family lies in the fact
that exponential distributions represent intervals between events
for Poisson processes on the real line and Freund distributions
can model bivariate processes with positive and negative
covariance. The Freund family of distributions becomes a
Riemannian 4-manifold with the Fisher information metric, and we
derive the induced $\alpha$-geometry, i.e., the $\alpha$-Ricci
curvature, the $\alpha$-scalar curvature etc. The case $\alpha=0$ recovers
the geometry of the metric or Levi Civita connection and we show that the
Freund manifold has a positive constant $0$-scalar curvature, so
geometrically it constitutes part of a sphere.

\subsection{Freund bivariate mixture exponential distributions }
Freund~\cite{freund} introduced a
bivariate exponential mixture distribution arising
from the following reliability considerations.  Suppose that an instrument has two
components $A$ and $B$ with lifetimes $X$ and $Y$ respectively having density
functions (when both components are in operation)\\
$f_{X}(x)=\alpha_{1}\,e^{-\alpha_{1} x}$;\\
$f_{Y}(y)=\alpha_{2}\,e^{-\alpha_{2} y}$ \\
for ($\alpha_{1},\alpha_{2}>0; x,y>0$).

Then $X$ and $Y$ are dependent in that a failure of
either component changes the parameter of the life distribution of
the other component.  Thus when $A$ fails, the parameter for $Y$
becomes $\beta_{2}$; when $B$ fails, the parameter for $X$ becomes
$\beta_{1}$.  There is no other dependence.  Hence the joint
density function of $X$ and $Y$ is :
\begin{eqnarray}
f(x,y)=\left\{\begin{array}{ll}
 \alpha_{1}\beta_{2}e^{-\beta_{2}y-(\alpha_{1}+\alpha_{2}-\beta_{2})x} &
\mbox{for}\, 0< x < y ,\\
\alpha_{2}\beta_{1}
e^{-\beta_{1}x-(\alpha_{1}+\alpha_{2}-\beta_{1})y} & \mbox{for}\,
0< y < x\, \end{array}\right. \label{freundbivariate}
\end{eqnarray}
where $\alpha_{i},\beta_{i} >0 \quad(i=1,2)$.

Provided that $\alpha_{1}+\alpha_{2}\neq \beta_{1},$ the marginal
density function of $X$ is
\begin{eqnarray}
f_{X}(x)&=& \left( \frac{\alpha_{2}}{
\alpha_{1}+\alpha_{2}-\beta_{1}} \right)\,\beta_{1}\,e^{-\beta_{1}
x} +\left(\frac{ \alpha_{1}
-\beta_{1}}{\alpha_{1}+\alpha_{2}-\beta_{1}} \right) \,
(\alpha_{1}+\alpha_{2})\,e^{-(\alpha_{1}+\alpha_{2})x}  \,,\,
x\geq0
\end{eqnarray}
and provided that $\alpha_{1}+\alpha_{2}\neq \beta_{2},$ The
marginal density function of $Y$ is
\begin{eqnarray}
f_{Y}(y)&=&\left(
\frac{\alpha_{1}}{\alpha_{1}+\alpha_{2}-\beta_{2}} \right)\,
\beta_{2}\, e^{-\beta_{2}y} + \left(
\frac{\alpha_{2}-\beta_{2}}{\alpha_{1}+\alpha_{2}-\beta_{2}}
\right)\, (\alpha_{1}+\alpha_{2})\, e^{-(\alpha_{1}+\alpha_{2})y}
\,,\, y\geq0
\end{eqnarray}
We can see that the marginal density functions are not exponential
but rather mixtures of exponential distributions if
$\alpha_{i}>\beta_{i}$ ; otherwise, they are weighted averages.
For this reason, this
system of distributions should be termed bivariate mixture
exponential distributions rather than simply bivariate exponential
distributions.  The marginal density functions $f_{X}(x)$ and
$f_{Y}(y)$ are exponential distributions only in the special case
$\alpha_{i}= \beta_{i} \quad (i=1,2)$.

Freund discussed the statistics of the special case when
$\alpha_{1}+\alpha_{2}= \beta_{1}= \beta_{2},$ and obtained the
joint density function as:

\begin{eqnarray}
f(x,y)=\left\{\begin{array}{ll}
 \alpha_{1}(\alpha_{1}+\alpha_{2})e^{-(\alpha_{1}+\alpha_{2})y} &
\mbox{for}\, 0< x < y ,\\
\alpha_{2}(\alpha_{1}+\alpha_{2})e^{-(\alpha_{1}+\alpha_{2})x} &
\mbox{for}\, 0< y < x\, \end{array}\right. \label{freundbivariatespecial}
\end{eqnarray}

with marginal density functions:
\begin{eqnarray}
f_{X}(x)&=& \left(\alpha_{1}+ \alpha_{2}(\alpha_{1}+\alpha_{2})x
\right)\, e^{-(\alpha_{1}+\alpha_{2})x} \,
x\geq0  \, ,\\
f_{Y}(y)&=&\left(\alpha_{2}+ \alpha_{1}(\alpha_{1}+\alpha_{2})y
\right)\, e^{-(\alpha_{1}+\alpha_{2})y} \, y\geq0
\end{eqnarray}

   The covariance and correlation coefficient of $X$ and $Y$ were derived by
   Freund, as follows:
\begin{eqnarray}
Cov(X,Y)&=& \frac{ {\beta_{1}}\,{\beta_{2}}-
{\alpha_{1}}\,{\alpha_{2}}}{{\beta_{1}}\,{\beta_{2}}\,{\left(
{\alpha_{1}} +
{\alpha_{2}} \right) }^2}\, ,\\
\rho(X,Y)&=&  \frac{ {\beta_{1}}\,{\beta_{2}}-
{\alpha_{1}}\,{\alpha_{2}} }
  {{\sqrt{{{\alpha_{2}}}^2 + 2\,{\alpha_{1}}\,{\alpha_{2}} + {{\beta_{1}}}^2}}\,{\sqrt{{{\alpha_{1}}}^2 + 2\,{\alpha_{1}}\,{\alpha_{2}} +
  {{\beta_{2}}}^2}}}
\end{eqnarray}
Note that\,$ -\frac{1}{3}< \rho(X,Y)<1 $.  The correlation
coefficient\,$ \rho(X,Y)\rightarrow 1 $ when  \,$ \beta_{1},\,
\beta_{2} \rightarrow \infty $, \, and $ \rho(X,Y)\rightarrow
-\frac{1}{3} $ when \,$\alpha_{1}=\alpha_{2}  $ and \,$
\beta_{1},\, \beta_{2} \rightarrow 0 $. In many applications,
$\beta_{i}>\alpha_{i} \quad (i=1,2)$ ( i.e., lifetime tends to be
shorter when the other component is out of action) ; in such cases
the correlation is positive.

\subsection{Fisher information metric}
 \begin{proposition}Let $F$ be the set of Freund bivariate mixture exponential
distributions, that is
 \begin{eqnarray}
 F= \{f
|f(x,y;\alpha_{1},\beta_{1},\alpha_{2},\beta_{2})=\left\{\begin{array}{ll}
 \alpha_{1}\beta_{2}e^{-\beta_{2}y-(\alpha_{1}+\alpha_{2}-\beta_{2})x} &
\mbox{{\rm for}}\, 0\leq x < y \\
\alpha_{2}\beta_{1}
e^{-\beta_{1}x-(\alpha_{1}+\alpha_{2}-\beta_{1})y} & \mbox{{\rm
for}}\, 0\leq y \leq x \end{array}\right.,\, \alpha_{i} ,
\beta_{i}>0\,(i,1,2)\} \label{freundmodel}
\end{eqnarray}
 Then we have :
\begin{enumerate}
\item Identifying $(\alpha_{1},\beta_{1},\alpha_{2},\beta_{2})$ as a local coordinate
system, $F$ can be regarded as a 4-manifold.
\item  $F$ becomes a Riemannian manifold with the Fisher information metric $
G=[g_{ij}]$ where
\begin{displaymath} g_{ij}=
\int_{0}^{\infty}\int_{0}^{\infty}\frac{\partial^ {2}\log
f(x,y)}{\partial x_{i}\partial x_{j}}\ f(x,y) \ dx\ dy
\end{displaymath}
and  $ (x_{1},\, x_{2},\, x_{3},\, x_{4})=(\alpha_{1}
,\,\beta_{1},\,\alpha_{2},\,
\beta_{2}) .$\\
is given by :
\begin{eqnarray}
  [g_{ij}]=\left[ \begin{array}{cccc}
   \frac{1}{{{\alpha_{1}}}^2 + {\alpha_{1}}\,{\alpha_{2}}} & 0 & 0 & 0 \\
   0 & \frac{{\alpha_{2}}}{{{\beta_{1}}}^2\,\left( {\alpha_{1}} + {\alpha_{2}} \right) } & 0 & 0 \\
   0 & 0 & \frac{1}{ {{\alpha_{2}}}^2+ {\alpha_{1}}\,{\alpha_{2}}} & 0 \\
   0 & 0 & 0 & \frac{{\alpha_{1}}}{{{\beta_{2}}}^2\, \left( {\alpha_{1}} + {\alpha_{2}} \right)}
 \end{array} \right] \label{freundmetric}
 \end{eqnarray}
\item The inverse $[g^{ij}]$ of  $[g_{ij}]$ is given by :
\begin{eqnarray}
   [g^{ij}]=\left[ \begin{array}{cccc}
   {{\alpha_{1}}}^2 + {\alpha_{1}}\,{\alpha_{2}} & 0 & 0 & 0 \\
    0 & \frac{{{\beta_{1}}}^2\,\left( {\alpha_{1}} + {\alpha_{2}} \right) }{{\alpha_{2}}} & 0 & 0 \\
    0 & 0 & {{\alpha_{2}}}^2+{\alpha_{1}}\,{\alpha_{2}} & 0 \\
     0 & 0 & 0 & \frac{{{\beta_{2}}}^2\, \left( {\alpha_{1}} +{\alpha_{2}} \right)}{{\alpha_{1}}}
  \end{array} \right]
  \end{eqnarray}
\end{enumerate}
$\hfill \Box$
\end{proposition}
\subsection{Geometry from the $\alpha$-connection}
We provide the $\alpha$-connection , and various
$\alpha$-curvature objects of the Freund manifold $F$: the
$\alpha$-curvature tensor, the $\alpha$-Ricci tensor, the
$\alpha$-scalar curvature, the $\alpha$-sectional curvatures and
the $\alpha$-mean curvatures.

\begin{enumerate}

\item $\alpha$-connection :

For each $ \alpha \in \mathbb{R}$, the {\em $\alpha$ (or
$\nabla^{(\alpha)}$)-connection} is the torsion-free affine
connection with components:

\begin{displaymath}
\Gamma_{ij,k}^{(\alpha)}= \int_{0}^{\infty} \int_{0}^{\infty}
\left( \frac{\partial^ {2}\log f}{\partial \xi^{i}\partial
\xi^{j}} \, \frac{\partial \log f}{\partial \xi^{k}}
+\frac{1-\alpha}{2}\,\frac{\partial \log f}{\partial
\xi^{i}}\,\frac{\partial \log f}{\partial \xi^{j}}\,\frac{\partial
\log f}{\partial \xi^{k}}\right) f \ dx \ dy
\end{displaymath}

\begin{proposition}The functions $\Gamma_{ij,k}^{(\alpha)}$ are given by :
\begin{eqnarray}
   \Gamma^{(\alpha)}_{11,1}&=&\frac{2\,\left(\alpha -1 \right) \,{{\alpha }_1} - \left( 1 + \alpha  \right) \,{{\alpha }_2}}
  {2\,{{{\alpha }_1}}^2\,{\left( {{\alpha }_1} + {{\alpha }_2} \right) }^2}
 \,,\nonumber \\
   \Gamma^{(\alpha)}_{11,3}&=&\frac{1 + \alpha }{2\,{{\alpha }_1}\,{\left( {{\alpha }_1} + {{\alpha }_2} \right) }^2}
 \,,\nonumber \\
   \Gamma^{(\alpha)}_{12,2}&=&\frac{\left( \alpha -1 \right)
   \,{{\alpha }_2}}{2\,{\left( {{\alpha }_1} + {{\alpha }_2} \right) }^2\,{{{\beta }_1}}^2}
   ,\nonumber \\
   \Gamma^{(\alpha)}_{13,3}&=&\frac{-1 + \alpha }{2\,{{\alpha }_2}\,{\left( {{\alpha }_1} + {{\alpha }_2} \right) }^2}
 \,,\nonumber \\
   \Gamma^{(\alpha)}_{14,4}&=&\frac{-\left( \left( \alpha -1 \right) \,{{\alpha }_2} \right) }
  {2\,{\left( {{\alpha }_1} + {{\alpha }_2} \right) }^2\,{{{\beta }_2}}^2}
   ,\nonumber \\
   \Gamma^{(\alpha)}_{22,2}&=&\frac{\left(  \alpha -1 \right)
    \,{{\alpha }_2}}{\left( {{\alpha }_1} + {{\alpha }_2} \right) \,{{{\beta }_1}}^3}
 \,,\nonumber \\
   \Gamma^{(\alpha)}_{22,3}&=&\frac{-\left( \left( 1 + \alpha  \right) \,{{\alpha }_1} \right) }
  {2\,{\left( {{\alpha }_1} + {{\alpha }_2} \right) }^2\,{{{\beta }_1}}^2}
 \,,\nonumber \\
   \Gamma^{(\alpha)}_{33,3}&=&\frac{-\left( \left( 1 + \alpha  \right) \,{{\alpha }_1} \right)  + 2\,\left( -1 + \alpha  \right) \,{{\alpha }_2}}
  {2\,{{{\alpha }_2}}^2\,{\left( {{\alpha }_1} + {{\alpha }_2} \right) }^2}
   ,\nonumber \\
   \Gamma^{(\alpha)}_{34,4}&=&\frac{\left(\alpha -1 \right)
   \,{{\alpha }_1}}{2\,{\left( {{\alpha }_1} + {{\alpha }_2} \right) }^2\,{{{\beta }_2}}^2}
   ,\nonumber \\
   \Gamma^{(\alpha)}_{44,4}&=&\frac{\left( \alpha -1 \right)
    \,{{\alpha }_1}}{\left( {{\alpha }_1} + {{\alpha }_2} \right) \,{{{\beta }_2}}^3}
 \end{eqnarray}
 while the other independent components are zero.\hfill $\Box$
\end{proposition}

We have an affine connection $\nabla^{(\alpha)}$ defined
 by :
 $$\langle  \nabla^{(\alpha)}_{\partial_i}\partial_j ,\partial_k \rangle=\Gamma_{ij,k}^{(\alpha)}\, ,$$
 So by solving the equations
$$  \Gamma_{ij,k}^{(\alpha)}= \sum_{h=1}^{4} g_{kh}\,\Gamma^{h(\alpha)}_{ij} \quad , (k=1,2,3,4).$$
we obtain the components of $\nabla^{(\alpha)}$ :

\begin{proposition}
The components $ \Gamma^{i(\alpha)}_{jk}$ of the
$\nabla^{(\alpha)}$-connections are given by :
\begin{eqnarray}
  \Gamma^{(\alpha)1}_{11}&=& \frac{1}{2}\,\left( - \frac{1 + \alpha }{{{\alpha }_1}}  +
   \frac{-1 + 3\,\alpha }{{{\alpha }_1} + {{\alpha }_2}}\right)\,,\nonumber\\
  \Gamma^{(\alpha)1}_{13}&=&\Gamma^{(\alpha)2}_{12}=\Gamma^{(\alpha)3}_{13}=\Gamma^{(\alpha)4}_{34}=
  \frac{-1 + \alpha }{2\,\left( {{\alpha }_1} + {{\alpha }_2} \right) } \, ,\nonumber\\
  \Gamma^{(\alpha)1}_{22}&=&-\Gamma^{(\alpha)3}_{22}= \frac{\left( 1 + \alpha  \right) \,{{\alpha }_1}\,{{\alpha }_2}}
  {2\,\left( {{\alpha }_1} + {{\alpha }_2} \right) \,{{{\beta }_1}}^2} \, ,\nonumber\\
  \Gamma^{(\alpha)1}_{33}&=& \Gamma^{(\alpha)2}_{23}= \frac{\left( 1 + \alpha  \right)
  \,{{\alpha }_1}}{2\,{{\alpha }_2}\,\left( {{\alpha }_1} + {{\alpha }_2} \right) } \,,\nonumber\\
  \Gamma^{(\alpha)1}_{44}&=&-\Gamma^{(\alpha)3}_{44}=\frac{-\left( \left( 1 + \alpha  \right) \,{{\alpha }_1}\,{{\alpha }_2} \right) }
  {2\,\left( {{\alpha }_1} + {{\alpha }_2} \right) \,{{{\beta }_2}}^2}  \, ,\nonumber\\
  \Gamma^{(\alpha)3}_{11}&=&\Gamma^{(\alpha)4}_{14}= \frac{\left( 1 + \alpha  \right)
  \,{{\alpha }_2}}{2\,{{\alpha }_1}\,\left( {{\alpha }_1} + {{\alpha }_2} \right) }  \,,\nonumber\\
  \Gamma^{(\alpha)3}_{33}&=&\frac{1}{2}\, \left(-\frac{1 + \alpha }{{{\alpha }_2}}
   + \frac{-1 + 3\,\alpha }{{{\alpha }_1} + {{\alpha }_2}}\right)  \,,\nonumber\\
 \Gamma^{(\alpha)4}_{44}&=& \frac{-1 + \alpha }{{{\beta }_2}} \,,
 \end{eqnarray}
 while the other independent components are zero.
 $\hfill \Box$
\end{proposition}

\item $\alpha$-Curvatures :

\begin{proposition}
The components $ R^{(\alpha)}_{ijkl}$ of the $\alpha$-curvature
tensor are given by:
\begin{eqnarray}
  R^{(\alpha)}_{1212}&= \frac{\left( {\alpha }^2 -1\right) \,{{{\alpha }_2}}^2}
  {4\,{{\alpha }_1}\,{\left( {{\alpha }_1} + {{\alpha }_2} \right) }^3\,{{{\beta }_1}}^2} \, ,\nonumber \\
   R^{(\alpha)}_{1223}&= \frac{\left( {\alpha }^2-1 \right)
   \,{{\alpha }_2}}{4\,{\left( {{\alpha }_1} + {{\alpha }_2} \right) }^3\,{{{\beta }_1}}^2}\,,\nonumber \\
    R^{(\alpha)}_{1414}&= \frac{\left( {\alpha }^2 -1\right)
    \,{{\alpha }_2}}{4\,{\left( {{\alpha }_1} + {{\alpha }_2} \right) }^3\,{{{\beta }_2}}^2} \,,\nonumber \\
     R^{(\alpha)}_{1434}&= \frac{-\left( {\alpha }^2-1 \right) \,{{\alpha }_1}}
  {4\,{\left( {{\alpha }_1} + {{\alpha }_2} \right) }^3\,{{{\beta }_2}}^2} \, ,\nonumber \\
      R^{(\alpha)}_{2323}&= \frac{\left( {\alpha }^2-1 \right)
      \,{{\alpha }_1}}{4\,{\left( {{\alpha }_1} + {{\alpha }_2} \right) }^3\,{{{\beta }_1}}^2}\, ,\nonumber \\
       R^{(\alpha)}_{2424}&=\frac{\left({\alpha }^2 -1\right) \,{{\alpha }_1}\,{{\alpha }_2}}
  {4\,{\left( {{\alpha }_1} + {{\alpha }_2} \right) }^2\,{{{\beta }_1}}^2\,{{{\beta }_2}}^2} \, ,\nonumber \\
       R^{(\alpha)}_{3434}&= \frac{\left( {\alpha }^2-1 \right) \,{{{\alpha }_1}}^2}
  {4\,{{\alpha }_2}\,{\left( {{\alpha }_1} + {{\alpha }_2} \right) }^3\,{{{\beta }_2}}^2} \,,
  \end{eqnarray}
while the other independent components are zero. $\hfill \Box$
\end{proposition}
Contracting $R^{(\alpha)}_{ijkl}$ with $g^{il}$ we obtain the
components $R^{(\alpha)}_{jk}$ of the $\alpha$-Ricci tensor.

\begin{proposition} The $\alpha$-Ricci tensor $ R^{(\alpha)}=[R^{(\alpha)}_{ij}]$ is given by :
\begin{eqnarray}
   R^{(\alpha)}=[R^{(\alpha)}_{jk}]=\left[ \begin{array}{cccc}
   \frac{-\left({\alpha }^2 -1\right) \,{{\alpha }_2}  }
    {2\,{{\alpha }_1}\,{\left( {{\alpha }_1} + {{\alpha }_2} \right) }^2} &
    0&
   \frac{{\alpha }^2-1}{2\,{\left( {{\alpha }_1} + {{\alpha }_2} \right) }^2}&
   0 \\
   0 & \frac{- \left( {\alpha }^2-1 \right) \,{{\alpha }_2} }
    {2\,\left( {{\alpha }_1} + {{\alpha }_2} \right) \,{{{\beta }_1}}^2}& 0 &
    0\\
  \frac{{\alpha }^2-1}{2\,{\left( {{\alpha }_1} + {{\alpha }_2} \right) }^2} &
  0&
   \frac{-\left({\alpha }^2 -1\right) \,{{\alpha }_1} }
    {2\,{{\alpha }_2}\,{\left( {{\alpha }_1} + {{\alpha }_2} \right) }^2}&
    0\\
   0 & 0 & 0 & \frac{-\left( {\alpha }^2-1 \right) \,{{\alpha }_1} }
    {2\,\left( {{\alpha }_1} + {{\alpha }_2} \right) \,{{{\beta }_2}}^2}
  \end{array} \right]
  \end{eqnarray}
The $\alpha$-eigenvalues and the $\alpha$-eigenvectors of the
$\alpha$-Ricci tensor are given by :
\begin{eqnarray}
  \left( {\alpha }^2-1 \right)  \left(\begin{array}{cccc}
  0 \\ \frac{1}{{\left( {{\alpha }_1} + {{\alpha }_2} \right) }^2}
  - \frac{1}{2\,{{\alpha }_1}\,{{\alpha }_2}}
       \\
  \frac{-{{\alpha }_2} }
   {2\,\left( {{\alpha }_1} + {{\alpha }_2} \right) \,{{{\beta
   }_1}}^2}\\
  \frac{-{{\alpha }_1}}
   {2\,\left( {{\alpha }_1} + {{\alpha }_2} \right) \,{{{\beta }_2}}^2}
   \end{array} \right)
  \end{eqnarray}

\begin{eqnarray}
   \left( \begin{array}{cccc}
 \frac{{{\alpha }_1}}{{{\alpha }_2}}& 0 & 1 & 0 \\
 - \frac{{{\alpha }_2}}{{{\alpha }_1}}  & 0 & 1& 0 \\
  0& 1& 0& 0\\
   0& 0& 0 & 1
  \end{array} \right)
  \end{eqnarray}
$\hfill \Box$
\end{proposition}

\begin{proposition} The manifold $F$ has a constant
$\alpha$-scalar curvature
\begin{eqnarray}
R^{(\alpha)}=\frac{-3\,\left({\alpha }^2-1 \right) }{2}
\end{eqnarray}
Note that the Freund manifold has a positive scalar curvature
$R^{(0)}=\frac{3}{2}$ when $\alpha=0$.  So geometrically it
constitutes part of sphere. $\hfill \Box$
\end{proposition}

 \begin{proposition} The $\alpha$-sectional curvatures $
\varrho^{(\alpha)}(\lambda,\mu)\,(\lambda,\mu=1,2,3,4)$ are given
by :
\begin{eqnarray}
  \varrho^{(\alpha)}(1,2)&=&\varrho^{(\alpha)}(1,4)=\frac{\left( 1 - {\alpha }^2 \right)
  \,{{\alpha }_2}}{4\,\left( {{\alpha }_1} + {{\alpha }_2} \right) } \, ,\nonumber \\
 \varrho^{(\alpha)}(1,3)&=& 0 ,\nonumber \\
  \varrho^{(\alpha)}(2,3)&=&\frac{\left( 1 - {\alpha }^2 \right) \,{{\alpha }_1}}
    {4\,\left( {{\alpha }_1} + {{\alpha }_2} \right) }  \,,\nonumber \\
  \varrho^{(\alpha)}(2,4)&=&  \frac{1 - {\alpha }^2}{4} \,,\nonumber \\
  \varrho^{(\alpha)}(3,4)&=& \varrho(2,3)\,.
 \end{eqnarray}
 $\hfill \Box$
 \end{proposition}

\begin{proposition}The $\alpha$-mean curvatures $
\varrho^{(\alpha)}(\lambda)\,(\lambda=1,2,3,4)$ are given by :
\begin{eqnarray}
 \varrho^{(\alpha)}(1)&=& \frac{\left( 1 - {\alpha }^2 \right)
 \,{{\alpha }_2}}{6\,\left( {{\alpha }_1} + {{\alpha }_2} \right) } \,  ,\nonumber \\
  \varrho^{(\alpha)}(2)&=&\varrho(4)= \frac{1 - {\alpha }^2}{6} \, ,\nonumber \\
  \varrho^{(\alpha)}(3)&=& \frac{\left( 1 - {\alpha }^2 \right) \,{{\alpha }_1}}
   {6\,\left( {{\alpha }_1} + {{\alpha }_2} \right) } \,.
\end{eqnarray}
$\hfill \Box$
 \end{proposition}

\end{enumerate}

\section{Submanifolds of the Freund 4-manifold}
We consider
five submanifolds $F_{i}\,(i=1,2,3,4)$ of the 4-manifold $F$ of
Freund bivariate exponential distributions
$f(x,y;\alpha_{1},\beta_{1},\alpha_{2},\beta_{2})$
(\ref{freundbivariate}), which includes the case of statistically
independent random variables. It includes also the special case of
an Absolutely Continuous Bivariate Exponential Distribution called
ACBED (or ACBVE) by Block and Basu (cf. Hutchinson and Lai
~\cite{hut}).  We use the coordinate system
$(\alpha_{1},\beta_{1},\alpha_{2},\beta_{2})$ for the submanifolds
$F_{i}\,(i\neq 3)$, and the coordinate system
$(\lambda_{1},\lambda_{12},\lambda_{2})$ for ACBED of the Block
and Basu case.

\subsection{Submanifold $ F_{1}\subset F $ : $\beta_{1}=\alpha_{1},$
$\beta_{2}=\alpha_{2}$}
 The distributions are of form :
 \begin{eqnarray}
f(x,y;\alpha_{1},\alpha_{2})=f_{1}(x;\alpha_{1})f_{2}(y;\alpha_{2})
 \end{eqnarray}
where $f_{i}$\, are the density functions of the one-dimensional
exponential  distributions with the parameters
\,$\alpha_{i}>0\,(i=1,2) $.  This is the case  for statistical
independence of $X$ and $Y$, so the space $ F_{1}$ is the direct
product of two Riemannian spaces
$\{f_{1}(x;\alpha_{1}):f_{1}(x;\alpha_{1})=\alpha_{1}
e^{-\alpha_{1}x},\,\alpha_{1}>0\}$ and
$\{f_{2}(y;\alpha_{2}):f_{2}(y;\alpha_{2})=-\alpha_{2}e^{-\alpha_{2}y},\,\alpha_{2}>0\}.$
\begin{proposition}
 The metric tensor
$[g_{ij}]$ is as follows :
\begin{eqnarray}
 [g_{ij}]=\left[ \begin{array}{cc}
 \frac{1}{\alpha_{1}^{2}} & 0 \\
 0 & \frac{1}{\alpha_{2}^{2}}
\end{array} \right]
\label{2freundmetric}
\end{eqnarray}
$\hfill \Box$
\end{proposition}
\begin{proposition}
The components of the $\alpha$-connection are
\begin{eqnarray}
\Gamma^{(\alpha)}_{11,1}&=&\frac{\alpha-1}{{\alpha_{1}}^3}
 \,,\nonumber \\
   \Gamma^{(\alpha)}_{22,2}&=&\frac{\alpha-1}{{\alpha_{2}}^3}
 \,,\nonumber \\
  \Gamma^{(\alpha)1}_{11}&=&\frac{\alpha-1}{{\alpha_{1}}}
  \,,\nonumber \\
  \Gamma^{(\alpha)2}_{22}&=&\frac{\alpha-1}{{\alpha_{2}}}\, ,
 \end{eqnarray}
 while the other components are zero.
 $\hfill \Box$
 \end{proposition}
\begin{proposition}
 The $\alpha$-curvature tensor, $\alpha$-Ricci tensor, and $\alpha$-scalar curvature of $F_1$ are
 zero.
 $\hfill \Box$
\end{proposition}

\subsection{Submanifold $ F_2\subset F $ : $\alpha_{1}=\alpha_{2},$
$\beta_{1}=\beta_{2}$ }
 The distributions are of form :
 \begin{eqnarray}
f(x,y;\alpha_{1},\beta_{1})=\left\{\begin{array}{ll}
 \alpha_{1}\beta_{1}\,
e^{-\beta_{1}y -(2\,\alpha_{1}-\beta_{1})x} &
\mbox{for}\, 0< x < y \\
 \alpha_{1}\beta_{1}\,e^{-\beta_{1}x-(2\,\alpha_{1}-\beta_{1})y} &
\mbox{for}\, 0< y < x\,
\end{array}\right.
 \end{eqnarray}
with parameters $\alpha_{1},\beta_{1}>0$.  The covariance,
correlation coefficient and marginal density functions, of $X$ and
$Y$ are given by :
\begin{eqnarray}
Cov(X,Y)&=& \frac{1}{4}\left(\frac{1}{{{{\alpha }_1}}^{2}} - \frac{1}{{{{\beta }_1}}^{2}}\right),\\
\rho(X,Y)&=& 1 - \frac{4\,{{{\alpha }_1}}^2}{3\,{{{\alpha }_1}}^2 + {{{\beta }_1}}^2},\\
f_{X}(x)&=& \left( \frac{\alpha_{1}}{2\, \alpha_{1}-\beta_{1}}
\right)\,\beta_{1}\,e^{-\beta_{1} x} +\left(\frac{ \alpha_{1}
-\beta_{1}}{2\,\alpha_{1}-\beta_{1}} \right) \,
(2\,\alpha_{1})\,e^{-2\,\alpha_{1}x}  \,,\, x\geq0 \,,\\
f_{Y}(y)&=&\left( \frac{\alpha_{1}}{2\,\alpha_{1}-\beta_{1}}
\right)\, \beta_{1}\, e^{-\beta_{1}y} + \left(
\frac{\alpha_{1}-\beta_{1}}{2\,\alpha_{1}-\beta_{1}} \right)\,
(2\,\alpha_{1})\, e^{-2\,\alpha_{1}y} \,,\, y\geq0 \,.
\end{eqnarray}
Note that $\rho(X,Y)=0$ when $\alpha_1=\beta_1.$

Note that $F_2$ forms as exponential family, with parameters
$(\alpha_1,\beta_1)$ and potential function
\begin{equation}
\psi= - \log(\alpha_1\, \beta_1)
  \end{equation}

 \begin{proposition}
The submanifold $F_2$ is an isometric isomorph of a the manifold $
F_1$.
\end{proposition}

{\bf Proof:}
 Since $\psi$ is a potential function, the Fisher
metric is given by the Hessian of $\psi$, that is,

\begin{equation}
 g_{ij}=\frac{\partial^{2}\psi}{\partial \theta^{i}\partial \theta^{j}}.
\end{equation}
Then, we have the Fisher metric (\ref{2freundmetric}) by a
straightforward calculation. \hfill $\Box$

\subsubsection{ Mutually dual foliations: } Since $
\nabla^{(1)}_{\partial_i}\partial_j=0 $, $(\alpha_1,\beta_1)$ is a
1-affine coordinate system, and the (-1)-affine coordinate system
is given by:
\begin{eqnarray}
  \eta_1 & = & -\frac{1}{\alpha_1}  , \nonumber \\
  \eta_2 & = & -\frac{1}{\beta_1} \,.
\end{eqnarray}
These coordinate have potential function:

\begin{equation}
  \lambda= \log(\alpha_1\, \beta_1)-2\, .
\end{equation}

So the coordinates $(\alpha_1,\beta_1)$ and
$(-\frac{1}{\alpha_1},-\frac{1}{\beta_1})$ are mutually dual with
respect to the Fisher metric, and the tetrad $ \{
F_2,g,\nabla^{(1)},\nabla^{(-1)} \}$ is dually flat space.
Therefore $F_2$ has dually orthogonal foliations.

For example: take $(\alpha_1,\eta_2)$ as a coordinate system for
$F_2$; then
\begin{eqnarray}
f(x,y;\alpha_{1},\eta_{2})=\left\{\begin{array}{ll}
 -\frac{\alpha_{1}}{\eta_2}\,
e^{\left( \frac{1}{\eta_2}\right)  y -\left( 2\,\alpha_1 +
\frac{1}{\eta_2} \right ) x} &
\mbox{for}\, 0< x < y \\
-\frac{\alpha_{1}}{\eta_2}\,e^{\left( \frac{1}{\eta_2} \right ) x
-\left( 2\,\alpha_1 + \frac{1}{\eta_2} \right ) y} & \mbox{for}\,
0< y < x\,
\end{array}\right.
 \end{eqnarray}

and the Fisher metric is :
\begin{eqnarray}
 [g_{ij}]=\left[ \begin{array}{cc}
 -\frac{1}{\alpha_{1}^{2}} & 0 \\
 0 & \frac{1}{(\eta_{2})^{2}}
\end{array} \right]
\end{eqnarray}

\subsubsection{Neighbourhoods of independence in $F_2$}
An important practical application of the Freund submanifold $F_2$
is the representation of a bivariate stochastic proces for which
the marginals are identical exponentials. The next result is
important because it provides topological neighbourhoods of  that
subspace $W$ in $F_2$ consisting of the bivariate processes that
have zero covariance: we obtain neighbourhoods of independence
for random (ie exponentially distributed) processes.

\begin{proposition}
Let $ \{ F_2,g,\nabla^{(1)},\nabla^{(-1)} \}$  be the manifold
$F_2$ with Fisher metric $g$ and exponential connection
$\nabla^{(1)}$.  Then $F_2$ can be realized in Euclidean ${\R}^3$
by the graph of a potential function, namely, $F_2$ can be
realized by the affine immersion $\{h,\xi\}$:
\[
  h: {\cal G} \rightarrow \R^3 :
  \left( \! \!
    \begin{array}{c}
       \alpha_1 \\ \beta_1
    \end{array} \! \! \right)
    \mapsto
    \left( \! \!
    \begin{array}{c}
       \alpha_1 \\ \beta_1 \\ \psi
    \end{array} \! \! \right), \quad
  \xi = \left( \! \! \begin{array}{c}
       0 \\ 0 \\ 1
    \end{array} \! \! \right).
\]

where $ \psi=- \log(\alpha_1\, \beta_1) $ and $\xi$ is the
transversal vector field along $ h.$

In $F_2$, the submanifold $W$ consisting of the independent case
$(\alpha_1= \beta_1) $ is represented by the curve :

$$(0,\infty)\rightarrow \R^3: (\alpha_1) \mapsto (\alpha_1 ,- \log(\alpha_1\, \beta_1)),
 \quad \xi=(0,0,1).$$
\end{proposition}

 This is illustrated in the graphic which shows $S$, an affine
 embedding of $F_2$ as a surface in $\R^3$, and $T$ an
 $\R^3$-tubular neighbourhood of $W$, the curve $\alpha_1= \beta_1$ in
 the surface.  This curve $W$ represents all bivariate distributions
 having identical exponential marginals and zero covariance; its
 tubular neighourhood $T$ represents departures from independence.

\begin{figure}
\begin{picture}(300,245)(0,0)
\put(-20,0){\includegraphics{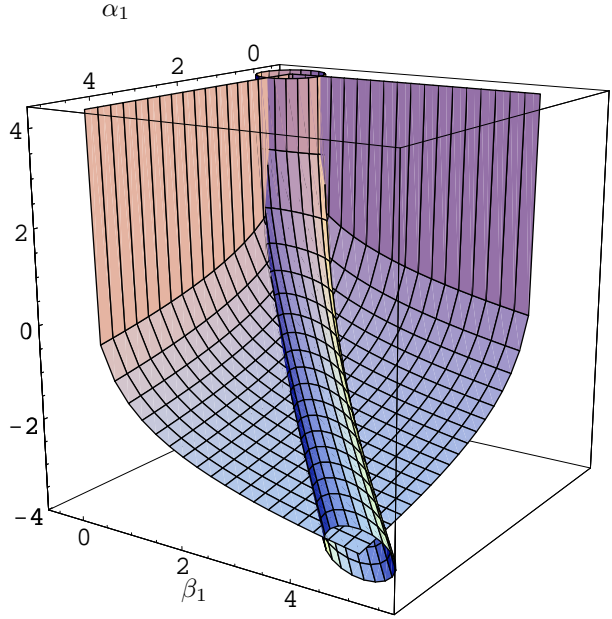}}
\put(140,10){$\beta_1$} \put(110,230){$\alpha_1$}
\end{picture}
\caption{{\em Affine immersion in natural coordinates
$(\alpha_1,\beta_1)$ as a surface in $\R^3$ for the Freund
submanifold $F_2;$ the tubular neighbourhood surrounds the
curve ($\alpha_1 = \beta_1$ in
the surface) consisting of all
bivariate distributions having identical exponential marginals and
zero covariance.}} \label{nhdF_2}
\end{figure}

\subsection{Submanifold $F_3 \subset F $: $\beta_{1}=\beta_{2}=\alpha_{1}+\alpha_{2} $}
The distributions are of form:
 \begin{eqnarray}
 f(x,y;\alpha_{1},\alpha_{2},\beta_{2})=\left\{\begin{array}{ll}
   \alpha_1\, \left(\alpha_1+\alpha_2 \right)\,e^{-(\alpha_1+\alpha_2)y}&
\mbox{for}\,\, 0< x < y \\
   \alpha_2\, \left(\alpha_1+\alpha_2 \right)\,e^{-(\alpha_1+\alpha_2)x}&
   \mbox{for}\,\,
0< y < x\, \end{array}\right.
\end{eqnarray}
with parameters $\alpha_{1}, \alpha_{2}>0$.  The covariance,
correlation coefficient and  marginal functions, of $X$ and $Y$
are given by :
\begin{eqnarray}
Cov(X,Y)&=&  \frac{{{{\alpha }_1}}^2 + {{\alpha }_1}\,{{\alpha
}_2} + {{{\alpha }_2}}^2}
  {{\left( {{\alpha }_1} + {{\alpha }_2} \right) }^4} ,\\
\rho(X,Y)&=&  \frac{{{{\alpha }_1}}^2 + {{\alpha }_1}\,{{\alpha
}_2} + {{{\alpha }_2}}^2}
  {\sqrt{2\,(\alpha_1+\alpha_2)^2-{\alpha_1}^2}\,\sqrt{2\,{\alpha_1}^2+4\,\alpha_1 \alpha_2+{\alpha_2}^2}} ,\\
f_{X}(x)&=& \left( {\alpha_2\, (\alpha_1+\alpha_2)x+\alpha_1 }\right) \,e^{-\left( {{\alpha }_1} + {{\alpha }_2} \right)x }  ,\,x\geq0\\
f_{Y}(y)&=& \left( {\alpha_1\, (\alpha_1+\alpha_2)y+\alpha_2
}\right) \,e^{-\left( {{\alpha }_1} + {{\alpha }_2} \right)y }
,\,y\geq0
 \end{eqnarray}
 Note that the correlation coefficient is positive.
 \begin{proposition}
 The metric tensor  $[g_{ij}]$ is given by :
\begin{eqnarray}
[g_{ij}]=\left[ \begin{array}{cc}
 \frac{{{\alpha }_2}+2\,{{\alpha }_1} }
    {{{\alpha }_1}\,{\left( {{\alpha }_1} + {{\alpha }_2} \right)
    }^2} & \frac{1}{{\left( {{\alpha }_1} + {{\alpha }_2} \right) }^{2}} \\
 \frac{1}{{\left( {{\alpha }_1} + {{\alpha }_2} \right) }^{2}} &
   \frac{{{\alpha }_1} + 2\,{{\alpha }_2}}
    {{{\alpha }_2}\,{\left( {{\alpha }_1} + {{\alpha }_2} \right) }^2}
\end{array} \right]
\end{eqnarray}
$\hfill \Box$
\end{proposition}
\begin{proposition}
The components of the $\alpha$-connection of $F_{3}$ are
\begin{eqnarray}
  \Gamma^{(\alpha)1}_{11}&=& \frac{-\left( \frac{1 + \alpha }{{{\alpha }_1}} \right)
   + \frac{-1 + 3\,\alpha }{{{\alpha }_1} + {{\alpha }_2}}}{2},\nonumber \\
  \Gamma^{(\alpha)1}_{12}&=& \frac{-1 + \alpha }{2\,\left( {{\alpha }_1} + {{\alpha }_2} \right) } ,\nonumber \\
  \Gamma^{(\alpha)1}_{22}&=& \frac{\left( 1 + \alpha  \right)
  \,{{\alpha }_1}}{2\,{{\alpha }_2}\,\left( {{\alpha }_1} + {{\alpha }_2} \right) } ,\nonumber \\
  \Gamma^{(\alpha)2}_{11}&=& \frac{\left( 1 + \alpha  \right)
  \,{{\alpha }_2}}{2\,{{\alpha }_1}\,\left( {{\alpha }_1} + {{\alpha }_2} \right) } ,\nonumber \\
   \Gamma^{(\alpha)2}_{22}&=& \frac{-\left( \frac{1 + \alpha }{{{\alpha }_2}} \right)  +
    \frac{-1 + 3\,\alpha }{{{\alpha }_1} + {{\alpha }_2}}}{2}  ,
 \end{eqnarray}
 while the other independent components are zero.
 $\hfill \Box$
\end{proposition}
\begin{proposition}
The $\alpha$-curvature tensor, $\alpha$-Ricci curvature, and
$\alpha$-scalar curvature of $F_3$ are zero.
\end{proposition}

\subsection{Submanifold $F_4 \subset F $: ACBED of Block and Basu}
Consider the distributions are form:
 \begin{eqnarray}
 f(x,y;\lambda_{1},\lambda_{12},\lambda_{2})=\left\{\begin{array}{ll}
 \frac{ {{\lambda }_1}\,{\lambda }\,\left( {{\lambda }_2} + {{\lambda }_{12}} \right) \,}{{{\lambda }_1} + {{\lambda }_2}}\,
    e^{- {{\lambda }_1}\,x  - \left( {{\lambda }_2}+ {{\lambda }_{12}} \right)\,y }\, &
\mbox{for}\,\, 0< x < y \\
 \frac{{{\lambda }_2}\,{\lambda } \,\left( {{\lambda }_1} + {{\lambda }_{12}} \right) }
    {{{\lambda }_1} + {{\lambda }_2}}\,
    e^{ -\left( {{\lambda }_1} + {{\lambda }_{12}} \right)\,x -{{\lambda }_2}\,y\,  }\, &
    \mbox{for}\,\,
0< y < x\, \end{array}\right.
\end{eqnarray}
where the parameters $\lambda_{1}, \lambda_{12}, \lambda_{2} $ are
positive, and $
\lambda= \lambda_{1}+ \lambda_{2}+ \lambda_{12}$.\\
This distribution was derived originally by omitting the singular
part of the Marshall and Olkin distribution (cf.~\cite{kotz}, page
[139]); Block and Basu called it the ACBED to emphasize that they
are the absolutely continuous bivariate exponential distributions.
Alternatively, it can be derived by Freund's method
(\ref{freundbivariate}), with
\begin{eqnarray}
 \alpha_{1}&=&\lambda_{1}+\frac{\lambda_{1}\,\lambda_{12}}{(\lambda_{1}+\lambda_{2})} ,\nonumber \\
 \beta_{1}&=&\lambda_{1}+\lambda_{12} ,\nonumber \\
 \alpha_{2}&=&\lambda_{2}+\frac{\lambda_{2}\,\lambda_{12}}{(\lambda_{1}+\lambda_{2})} ,\nonumber \\
 \beta_{2}&=&\lambda_{2}+\lambda_{12}. ,\nonumber
  \end{eqnarray}
 By substitutions we obtained the covariance, correlation coefficient
 and  marginal functions:
\begin{eqnarray}
 Cov(X,Y)&=&  \frac{{\left( {{\lambda }_1} + {{\lambda }_2} \right) }^2 \left( {{\lambda }_1} + {{\lambda }_{12}} \right)
   \left( {{\lambda }_2} + {{\lambda }_{12}} \right)-{{\lambda }}^2\,{{\lambda
   }_1}\,{{\lambda }_2}}
     {{{\lambda }}^2\,{\left( {{\lambda }_1} + {{\lambda }_2} \right) }^2\,\left( {{\lambda }_1} + {{\lambda }_{12}} \right) \,
       \left( {{\lambda }_2} + {{\lambda }_{12}} \right) } ,\\
       \rho(X,Y)&=&\frac{ {\left( {{\lambda }_1} + {{\lambda }_2} \right) }^2 \left( {{\lambda }_1} +
        {{\lambda }_{12}} \right)
   \left( {{\lambda }_2} + {{\lambda }_{12}} \right)-{{\lambda }}^2\,{{\lambda
   }_1}\,{{\lambda }_2}}
      {\sqrt{ \prod_{i=1,\, j\neq i }^{ 2} \left( {( {{\lambda }_1} + {{\lambda }_2}
      )}^2
     { ( {{\lambda }_i} + {{\lambda }_{12}} ) }^2 +
       {{\lambda }_j}{\lambda }^2 \left( {{\lambda }_j} + 2{{\lambda }_i} \right)
       \right)
       } },\\
 f_{X}(x)&=& \left(\frac{-{\lambda_{12} }}{ {{\lambda }_1} + {{\lambda }_{2}}
} \right)\,{\lambda }\,e^{-{\lambda }\,x }+ \left(\frac{{\lambda
}}{ {{\lambda }_1} + {{\lambda }_2} } \right)\,\left( {{\lambda
}_1} + {{\lambda }_{12}}
 \right)\,e^{-\left( {{\lambda }_1} + {{\lambda }_{12}} \right)\,x }
 ,\,x\geq0  \label{xmargial} \\
 f_{Y}(y)&=& \left(\frac{-{\lambda_{12} }}{ {{\lambda }_1} + {{\lambda }_{2}}} \right)\,{\lambda }\,e^{-{\lambda }\,y }+
  \left(\frac{{\lambda }}{{{\lambda }_1} + {{\lambda }_{2}}} \right) \,\left( {{\lambda }_2} + {{\lambda }_{12}}
 \right)\,e^{-\left( {{\lambda }_2} + {{\lambda }_{12}}
 \right)\,y }  ,\,y\geq0  \label{ymarginal}
 \end{eqnarray}
 Note that the correlation coefficient is positive, and the marginal functions are a negative mixture of two exponentials.
 \begin{proposition}
 The metric tensor  $[g_{ij}]$ using the coordinate system $ ({\lambda }_1,{\lambda }_{12},{\lambda }_2)$ is as follows:
\begin{eqnarray}
[g_{ij}]=\left[\begin{array}{ccc}
    \frac{{{\lambda }_2}\,\left( \frac{1}{{{\lambda }_1}} +
         \frac{{{\lambda }_1} + {{\lambda }_2}}{{\left( {{\lambda }_1} + {{\lambda }_{12}} \right) }^2} \right) }{{\left(
            {{\lambda }_1} + {{\lambda }_2} \right) }^2}+\frac{1}{{{\lambda } }^{2}}  &
   \frac{{{\lambda }_2}}{\left( {{\lambda }_1} + {{\lambda }_2} \right) \,
       {\left( {{\lambda }_1} + {{\lambda }_{12}} \right) }^2} +
    \frac{1}{{{\lambda } }^{2}} &
   \frac{-1}{{\left( {{\lambda }_1} + {{\lambda }_2} \right) }^{2}} +
    \frac{1}{{{\lambda } }^{2}} \\
   \frac{{{\lambda }_2}}
     {\left( {{\lambda }_1} + {{\lambda }_2} \right) \,{\left( {{\lambda }_1} + {{\lambda }_{12}} \right) }^2} +
    \frac{1}{{{\lambda } }^{2}}&
    \frac{\frac{{{\lambda }_2}}{{\left( {{\lambda }_1} + {{\lambda }_{12}} \right) }^2} +
       \frac{{{\lambda }_1}}{{\left( {{\lambda }_2} + {{\lambda }_{12}} \right) }^2}}{{{\lambda }_1} + {{\lambda}_2}}+\frac{1}{{{\lambda } }^{2}} &
   \frac{{{\lambda }_1}}{\left( {{\lambda }_1} + {{\lambda }_2} \right) \,
       {\left( {{\lambda }_2} + {{\lambda }_{12}} \right) }^2} +
    \frac{1}{{{\lambda } }^{2}}  \\
   \frac{-1}{{\left( {{\lambda }_1} + {{\lambda }_2} \right) }^{2} }+
    \frac{1}{{{\lambda } }^{2}}&
   \frac{{{\lambda }_1}}{\left( {{\lambda }_1} + {{\lambda }_2} \right) \,
       {\left( {{\lambda }_2} + {{\lambda }_{12}} \right) }^2} +
    \frac{1}{{{\lambda } }^{2}}&
    \frac{{{\lambda }_1}\,\left( \frac{1}{{{\lambda }_2}} +
         \frac{{{\lambda }_1} + {{\lambda }_2}}{{\left( {{\lambda }_2} + {{\lambda }_{12}} \right) }^2} \right) }{{\left(
            {{\lambda }_1} + {{\lambda }_2} \right) }^2}+\frac{1}{{{\lambda }
            }^{2}}
\end{array} \right]
\end{eqnarray}
$\hfill \Box$
\end{proposition}

The Christoffel symbols, curvature tensor, Ricci tensor, scalar
curvature, sectional curvatures and the mean curvatures were
computed but are not listed
because they are somewhat cumbersome.\\

In the case when $ {\lambda }_1={\lambda }_2 $, this family of
distributions becomes
\begin{eqnarray}
 f(x,y;\lambda_{1},\lambda_{12})=\left\{\begin{array}{ll}
 \frac{ (2\,\lambda_{1}+\lambda_{12})\,\left( {{\lambda }_1} +
  {{\lambda }_{12}} \right) }{2}\,
    e^{- {{\lambda }_1}\,x  - \left( {{\lambda }_1}+ {{\lambda }_{12}} \right)\,y }\, &
\mbox{for}\,\, 0< x < y \\
 \frac{(2\,\lambda_{1}+\lambda_{12})\,\left( {{\lambda }_1} +
  {{\lambda }_{12}} \right) }{2}\,e^{- {{\lambda }_1}\,y  - \left( {{\lambda }_1}+
  {{\lambda }_{12}} \right)\,x }\, & \mbox{for}\,\,
0< y < x\, \end{array}\right.
\end{eqnarray}
which is an exponential family with natural parameters
$(\theta_1,\,\theta_{2})= ({\lambda }_1,\,{\lambda }_{12}) $ and
potential function $\psi(\theta)= \log(2)-\log({\lambda
}_1+{\lambda }_{12})-\log(2\,{\lambda }_1+{\lambda }_{12})$, note
that from equations (\ref{xmargial}, \ref{ymarginal}), this family
of distributions has an identical marginal density functions.

So it would be easy to derive the $ \alpha$
-geometry, for example:\\

The metric tensor $[g_{ij}]$ is as follows:
\begin{eqnarray}
 [g_{ij}]=\left[ \begin{array}{cc}
 \frac{1}{{\left( {{\lambda }_1} + {{\lambda }_{12}} \right) }^{2}} +
    \frac{4}{{\left( 2\,{{\lambda }_1} + {{\lambda }_{12}} \right)
    }^2}&
   \frac{1}{{\left( {{\lambda }_1} + {{\lambda }_{12}} \right) }^{2}} +
    \frac{2}{{\left( 2\,{{\lambda }_1} + {{\lambda }_{12}} \right) }^2}\\
    \frac{1}{ {\left( {{\lambda }_1} + {{\lambda }_{12}} \right) }^{2}} +
    \frac{2}{{\left( 2\,{{\lambda }_1} + {{\lambda }_{12}} \right)
    }^2}&
   \frac{1}{{\left( {{\lambda }_1} + {{\lambda }_{12}} \right) }^{2}} +
  \frac{1}{ {\left( 2\,{{\lambda }_1} + {{\lambda }_{12}} \right)
  }^{2}} \end{array} \right]
\end{eqnarray}

By direct calculation the functions $\Gamma_{ij,k}^{(\alpha)}=
\frac{1-\alpha}{2}\,
{\partial_{i}\,\partial_{j}\,\partial_{k}\psi(\theta)};$ are given
by:
\begin{eqnarray}
 \Gamma^{(\alpha)}_{11,1}&=&{\left( 1 - \alpha  \right) \,\left( \frac{-1}
 {{\left( {{\lambda }_1} + {{\lambda }_{12}} \right) }^3} -
         \frac{8}{{\left( 2\,{{\lambda }_1} + {{\lambda }_{12}} \right) }^3} \right) }
 \,,\nonumber \\
   \Gamma^{(\alpha)}_{11,2}&=&{\left( 1 - \alpha  \right) \,
   \left( \frac{-1}{{\left( {{\lambda }_1} + {{\lambda }_{12}} \right) }^3} -
         \frac{4}{{\left( 2\,{{\lambda }_1} + {{\lambda }_{12}} \right) }^3} \right) }
 \,,\nonumber \\
   \Gamma^{(\alpha)}_{12,2}&=&{\left( 1 - \alpha  \right) \,\left( \frac{-1}{{\left( {{\lambda }_1} + {{\lambda }_{12}} \right) }^3} -
         \frac{2}{{\left( 2\,{{\lambda }_1} + {{\lambda }_{12}} \right) }^3} \right) }
         \,,\nonumber \\
   \Gamma^{(\alpha)}_{22,2}&=&{\left( 1 - \alpha  \right) \,\left( \frac{-1}{{\left( {{\lambda }_1} + {{\lambda }_{12}} \right) }^3} -
         \frac{1}{{\left( 2\,{{\lambda }_1} + {{\lambda }_{12}} \right) }^3} \right)
         }\,.
  \end{eqnarray}

 By solving the equations
$$  \Gamma_{ij,k}^{(\alpha)}= \sum_{h=1}^{3} g_{kh}\,\Gamma^{h(\alpha)}_{ij} \quad , (k=1,2)$$
we obtain the components of $\nabla^{(\alpha)}$ as follows:
\begin{eqnarray*}
 \Gamma^{(\alpha)1}=[\Gamma^{(\alpha)1}_{ij}]=\left[ \begin{array}{cc}
 \frac{1 - \alpha }{{{\lambda }_1} + {{\lambda }_{12}}} +
    \frac{4\,\left( \alpha-1  \right) }{2\,{{\lambda }_1} + {{\lambda
    }_{12}}}&
   \frac{\left(  \alpha-1  \right) \,{{\lambda }_{12}}}
    {\left( {{\lambda }_1} + {{\lambda }_{12}} \right) \,\left( 2\,{{\lambda }_1} + {{\lambda }_{12}} \right) }\\
     \frac{\left( \alpha -1 \right) \,{{\lambda }_{12}}}
    {\left( {{\lambda }_1} + {{\lambda }_{12}} \right) \,\left( 2\,{{\lambda }_1} + {{\lambda }_{12}} \right)
    }&
   \frac{-\left( \alpha -1  \right) \,{{\lambda }_1}}
      {\left( {{\lambda }_1} + {{\lambda }_{12}} \right) \,\left( 2\,{{\lambda }_1} +
       {{\lambda }_{12}} \right) }
\end{array} \right]   \, ,\nonumber
\end{eqnarray*}
\begin{eqnarray}
  \Gamma^{(\alpha)2}=[\Gamma^{(\alpha)2}_{ij}]=\left[ \begin{array}{cc}
 \frac{-2\,\left(\alpha-1  \right) \,{{\lambda }_{12}}}
    {\left( {{\lambda }_1} + {{\lambda }_{12}} \right) \,\left( 2\,{{\lambda }_1} + {{\lambda }_{12}} \right)
    }&
   \frac{2\,\left( \alpha -1 \right) \,{{\lambda }_1}}
    {\left( {{\lambda }_1} + {{\lambda }_{12}} \right) \,\left( 2\,{{\lambda }_1} + {{\lambda }_{12}} \right) }\\
     \frac{2\,\left( \alpha -1 \right) \,{{\lambda }_1}}
    {\left( {{\lambda }_1} + {{\lambda }_{12}} \right) \,\left( 2\,{{\lambda }_1} + {{\lambda }_{12}} \right)
    }&
   \frac{2\,\left( \alpha -1 \right) }{{{\lambda }_1} + {{\lambda }_{12}}} +
    \frac{1 - \alpha }{2\,{{\lambda }_1} + {{\lambda }_{12}}}
\end{array} \right]   \, .\nonumber
\end{eqnarray}

In this case, the $\alpha$-curvature tensor, $\alpha$-Ricci
curvature, and $\alpha$-scalar curvature are zero.\\

In addition, since the coordinates $({{\lambda }_1},{\lambda
}_{12})$ is a 1-affine coordinate system, then (-1)-affine
coordinate system is
$$(\eta_1,\eta_2)= (-\frac{1}{{{\lambda }_1} + {{\lambda }_{12}}}
- \frac{1}{{{\lambda }_1} + 2\,{{\lambda }_{12}}},
  - \frac{1}{{{\lambda }_1} + {{\lambda }_{12}}}
   - \frac{1}{2\,{{\lambda }_1} + {{\lambda }_{12}}} )$$  with potential function
   $$ \lambda=-2 -\log (2)+ \log (2\,{{\lambda }_1} + {{\lambda }_{12}}) + \log
({{\lambda }_1} + {{\lambda }_{12}}) .$$

\section{Freund bivariate mixture log-exponential distributions}
In this section we introduce a Freund bivariate mixture
log-exponential distribution which has mixture log-exponential
marginal functions, and discus their properties.

The Freund bivariate mixture log-exponential distributions arise
from the Freund distributions
(\ref{freundbivariate})for the non-negative random variables $x=
log \frac{1}{n}$ and
$y=log\frac{1}{m}$, or equivalently, $n=e^{-x}$ and $m=e^{-y}$.\\
So the Freund log-exponential distributions are
given by:
\begin{eqnarray}
g(n,m)=\left\{\begin{array}{ll}
 {\alpha_{1}\,\beta_{2}} \,{ m^{(\beta_{2}-1)}\,n^{(\alpha_{1}+\alpha_{2}-\beta_{2}-1)}} &
\mbox{for}\, 0<m < n<1 ,\\
{\alpha_{2}\,\beta_{1}}\,{
n^{(\beta_{1}-1)}\,m^{(\alpha_{1}+\alpha_{2}-\beta_{1}-1)}} &
\mbox{for}\, 0<n < m< 1\, \end{array}\right.
\end{eqnarray}
where $\alpha_{i}, \beta_{i} >0\quad (i=1,2)$.  The covariance,
and marginal density functions, of $n$ and $m$ are given by:
\begin{eqnarray}
Cov(n,m)&=&  \frac{{{\alpha }_2}\,\left( -\left( {{\alpha
}_1}\,\left( 2 + {{\alpha }_1} + {{\alpha }_2} \right)  \right)  +
       {{\beta }_1} \right)  + \left( {{\alpha }_1} +
       \left( {{\alpha }_1} + {{\alpha }_2} \right) \,{{\beta }_1} \right) \,{{\beta }_2}}{{\left( 1 +
        {{\alpha }_1} + {{\alpha }_2} \right) }^2\,\left( 2 + {{\alpha }_1} + {{\alpha }_2} \right) \,
    \left( 1 + {{\beta }_1} \right) \,\left( 1 + {{\beta }_2} \right) } \,,\\
g_{N}(n)&=& \left( \frac{\alpha_{2}}{
\alpha_{1}+\alpha_{2}-\beta_{1}} \right)\,\beta_{1}n^{\beta_{1}
-1} +\left(\frac{ \alpha_{1}
-\beta_{1}}{\alpha_{1}+\alpha_{2}-\beta_{1}} \right)
\,(\alpha_{1}+\alpha_{2})n^{(\alpha_{1}+\alpha_{2})-1}  \,,\\
g_{M}(m)&=&\left(
\frac{\alpha_{1}}{\alpha_{1}+\alpha_{2}-\beta_{2}} \right)
\,\beta_{2} m^{\beta_{2}-1} + \left(
\frac{\alpha_{2}-\beta_{2}}{\alpha_{1}+\alpha_{2}-\beta_{2}}
\right) \,(\alpha_{1}+\alpha_{2}) m^{(\alpha_{1}+\alpha_{2})-1}
\,.
\end{eqnarray}
Note that the marginal functions are mixture log-exponential
distributions. Directly from the definition of the Fisher metric we deduce:
\begin{proposition}
The family of Freund bivariate mixture log-exponential
distributions for random variables $n,m$ determines a Riemannian
4-manifold which is an isometric isomorph of the Freund 4-manifold .
\end{proposition}

\section{Concluding remarks}
We have derived the information geometry of the 4-manifold of
Freund bivariate mixture exponential distributions, which admits
positive and negative covariance.  The curvature objects are
derived and so also are those on four submanifolds, including the case of
statistically independent random variables, and the special case
ACBED of Block and Basu. We use one submanifold to
provide examples of neighbourhoods of the independent case for
bivariate distributions having identical exponential marginals.
Thus, since exponential distributions complement Poisson point
processes, we obtain a means to discuss the neighbourhood of
independence for random processes in general. The Freund manifold has a constant
$0$-scalar curvature, so geometrically it constitutes part of a
sphere.

The authors used {\em Mathematica} to perform analytic calculations~\cite{MJ03},
and can make available working notebooks for others to use.

\textbf{\textbf{Acknowledgment:}}The authors wish to thank the
referees for suggestiog improvements and the Libyan Ministry of
Education for a scholarship for Arwini.

\end{document}